\documentclass{article}

     \usepackage[final]{style/tackling_climate_workshop_style}

\usepackage[utf8]{inputenc} % allow utf-8 input
\usepackage[T1]{fontenc}    % use 8-bit T1 fonts
\usepackage{hyperref}       % hyperlinks
\usepackage{url}            % simple URL typesetting
\usepackage{booktabs}       % professional-quality tables
\usepackage{amsfonts}       % blackboard math symbols
\usepackage{nicefrac}       % compact symbols for 1/2, etc.
\usepackage{microtype}      % microtypography
\usepackage{multirow}

\usepackage{float}
\usepackage{stfloats}
\usepackage{lipsum}
\usepackage{graphicx,trimclip}
\usepackage{float}
\usepackage[table]{xcolor}
\usepackage{enumitem}

\usepackage{todonotes} % to be removed in final version

\title{Optimizing Energy Efficiency and Grid Stability via Public EV Charging Flexibility}

\author{%
  Marek Miltner\thanks{Alternative email: marek.miltner@cvut.cz} \\
  FEE CTU\\
  Prague, Czechia;\\
  CEE, Stanford University\\
    Stanford, USA\\
  \texttt{marek.miltner@stanford.edu} \\
  \And
  Artem Bryksa \\
  FEE CTU \\
  Prague, Czechia\\
  \texttt{bryksart@fel.cvut.cz} \\
  \AND
  Ondřej Štogl \\
  FEE CTU \\
  Prague, Czechia\\
  \texttt{stoglond@fel.cvut.cz} \\
  \And
  Jakub Zíka \\
  FEE CTU \\
  Prague, Czechia\\
  \texttt{zikajak3@fel.cvut.cz} \\
    \And
  Daniel Vašata \\
  FIT CTU \\
  Prague, Czechia\\
  \texttt{Daniel.Vasata@fit.cvut.cz} \\
    \And
  Magda Friedjungová \\
  FIT CTU \\
  Prague, Czechia\\
  \texttt{Magda.Friedjungova@fit.cvut.cz} \\
  \And
  Ram Rajagopal \\
  CEE, Stanford University\\
    Stanford, USA\\
  \texttt{ramr@stanford.edu} \\
    \And
  Oldřich Starý \\
  FEE CTU \\
  Prague, Czechia\\
  \texttt{staryo@fel.cvut.cz} \\
}

\begin{document}

\maketitle

\begin{abstract}
 This study evaluates the potential of electric vehicle (EV) charging flexibility to enhance both energy efficiency and power grid stability. Using real-world data from public charging stations in Prague, we analyze individual and aggregated charging sessions to explore how optimizing charging times can reduce energy waste, minimize grid imbalances, and support the integration of renewable energy. By aligning EV charging with periods of lower grid demand and higher renewable generation, we demonstrate a significant improvement in energy efficiency, reducing the need for costly system support and ancillary services. Our findings suggest that cooperation between power distributors and transmission system operators can unlock new opportunities for maintaining grid stability while promoting sustainable energy use in an increasingly uncertain energy landscape.

\end{abstract}

\section{Introduction and motivation}

As the share of renewable energy sources (RES) grows within global power systems in line with the green transformation of the industry \cite{pedro_opportunities_2023}, maintaining grid stability while optimizing energy efficiency has become a critical challenge for transmission system operators (TSOs) \cite{basit2020limitations}. The intermittent nature of RES generation, particularly from sources like wind and solar, introduces significant variability in electricity supply, often resulting in imbalances between supply and demand. Addressing these imbalances promptly and efficiently is crucial for avoiding unnecessary energy waste and ensuring that energy usage aligns with grid needs \cite{SAHA2023108701,SINSEL20202271}. 

TSOs play a vital role in ensuring both the reliable and energy-efficient operation of the transmission system. This responsibility includes the safe and efficient transmission of electricity from generators to distributors, real-time balancing of the grid, and overseeing essential processes such as frequency restoration and voltage control. Achieving these objectives while minimizing energy waste is increasingly difficult as RES penetration grows. Traditional storage solutions and flexibility resources are becoming inadequate, leading to rising costs in system services (SyS) and ancillary services (AnS) designed to maintain grid stability \cite{majer2010system}.

In this context, the need for innovative approaches to balance the grid efficiently, while minimizing energy losses and enhancing the utilization of renewable energy, is paramount. The ability to dynamically manage both generation and demand through advanced technologies can significantly contribute to overall energy efficiency and system reliability \cite{saele_understanding_2023}.

% \newpage
\section{EV Charging Flexibility as a Cost-Effective Approach to Energy Efficiency and Grid Stability}

One of the key challenges in achieving energy efficiency within power systems is the lack of sufficient storage capacity to balance the grid as the share of intermittent renewable energy sources (RES) increases. This challenge stems from the inherent variability of RES, such as wind and solar power, which do not generate a constant supply of electricity. As RES penetration grows, balancing supply and demand in real-time becomes increasingly complex and energy waste becomes more pronounced, especially during periods of surplus generation or when demand peaks cannot be met efficiently. Without flexible storage or demand-side solutions, a significant portion of the generated renewable energy is at risk of being curtailed, leading to inefficiencies in both energy use and economic cost.

This article focuses on the untapped flexibility potential of existing battery capacity—specifically electric vehicles (EVs) connected to public charging stations—as a solution that requires no major central investment. Unlike large-scale grid storage, which demands substantial infrastructure, EV batteries offer a decentralized, scalable resource already in place. By leveraging the flexibility of EV charging, we can optimize energy usage by matching charging demand to times when renewable energy is abundant, thus reducing both grid imbalances and energy waste \cite{tucki_development_2019, shang_fedpt-v2g_2024}. This approach not only improves the overall efficiency of energy use but also helps lower operational costs for transmission system operators (TSOs) by reducing the need for expensive ancillary services and system support mechanisms.

Our assessment is based on real-world data from public charging points in Prague, Czech Republic. By analyzing detailed information such as connection and disconnection times, the amount of energy transferred, and average charging power, this study estimates the potential of EVs to serve as dynamic, flexible resources for grid stabilization. These data provide insights into how individual EV charging sessions can be adapted to better align with grid conditions, maximizing the use of renewable energy and minimizing the need for non-renewable peaking power plants. Additionally, by cross-referencing these charging data with grid supply-demand imbalances and flexibility requests issued by ČEPS, the Czech TSO, we are able to evaluate how effectively EVs can be used to support grid efficiency. This integration of EVs as distributed storage assets contributes to a more resilient and energy-efficient power system.

\subsection{Individual Flexibility}

The core concept of EV charging flexibility begins at the level of individual charging sessions. Each session has a set start and end time, based on when the vehicle is connected to the charger. However, many vehicles do not charge for the entire duration of the session, as they often reach full charge early, leaving a window of idle connection time that can be exploited for energy efficiency gains \cite{develder2016quantifying}. For example, a vehicle may be plugged in overnight, but it only requires a few hours of actual charging time. The remaining time represents a significant opportunity to shift the load to periods when renewable energy generation is higher or when grid demand is lower.

Currently, EVs typically charge at the maximum permissible rate from the start of the session until full charge is achieved. This inflexible charging pattern often coincides with peak grid demand periods, exacerbating grid stress and leading to higher energy costs and inefficiencies. However, from an energy efficiency perspective, this is not the most optimal approach. If we allow flexibility in the timing of the charging within the session window, we can create multiple alternative charging load profiles that reduce peak demand on the grid while ensuring that the vehicle still reaches full charge by the time it is needed.

\begin{figure}[H]
  \centering

  \includegraphics[width=1\textwidth]{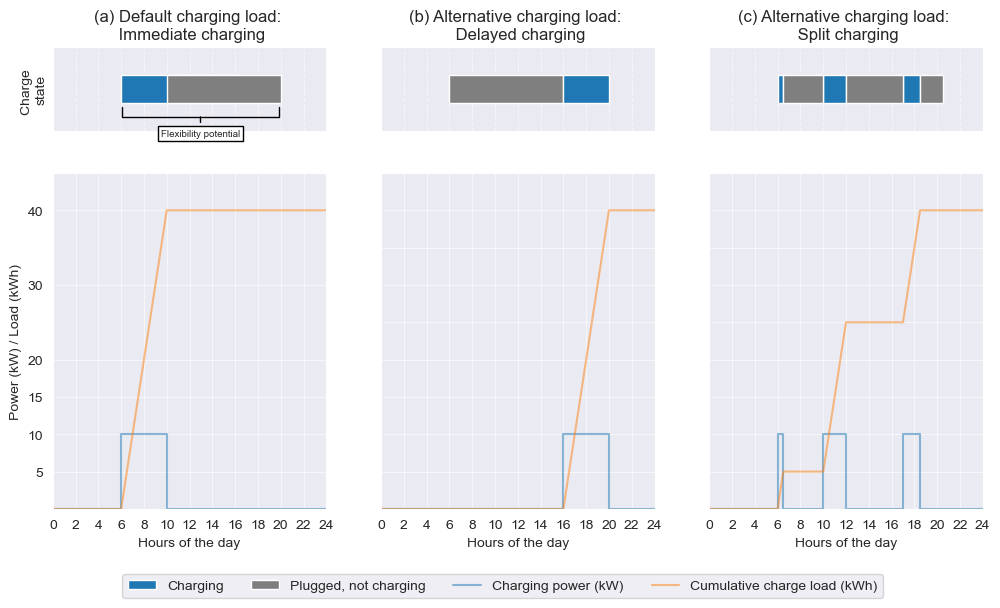}

  \caption{Individual EV charging load flexibility demonstration with three equivalent outcomes from the user’s point of view, each with different impacts on grid efficiency. Flexibility allows charging to align with energy-efficient periods of low grid demand or high renewable generation.}
  \label{gantt-fig}
\end{figure}

For example, by delaying the charging to start later in the session or by modulating the charging power based on real-time grid conditions, EVs can shift their consumption to periods when electricity is cheaper, greener, or more abundant. Figure \ref{gantt-fig} shows various potential load profiles, all of which provide the same final state of charge to the EV owner but differ in their impact on the grid. By shifting load times, we can avoid the use of inefficient peaking plants, which are typically fossil-fuel-based and expensive to operate, and improve grid efficiency without affecting the user experience. This flexibility is key to integrating higher levels of renewable energy into the grid without sacrificing stability.

\subsection{Potential of Aggregated EV Charging Flexibility}
\label{potential-flex}

The potential benefits of aggregated EV charging flexibility are even greater when multiple vehicles are considered together. By coordinating the charging behavior of many vehicles, system operators can smooth out overall demand and better align energy consumption with periods of renewable energy availability. Aggregating EVs as a flexible resource allows grid operators to treat them as a "virtual power plant" that can absorb excess renewable energy or reduce demand during peak hours, thereby significantly enhancing both energy efficiency and grid stability.

\begin{figure}[H]
  \centering

  \includegraphics[width=1\textwidth]{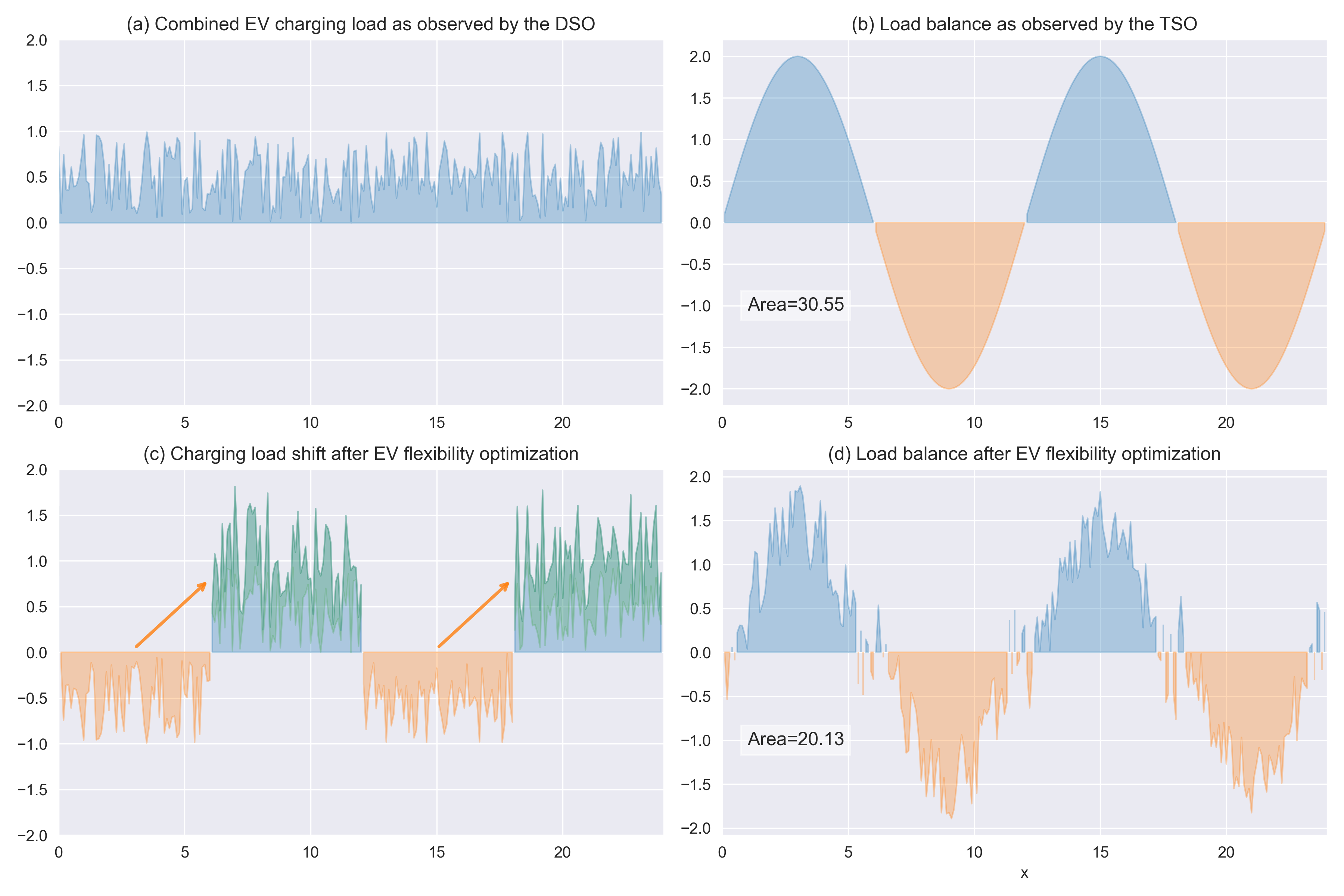}

  \caption{Subplots (a) and (b) show observed charging load and system imbalance, while subplots (c) and (d) demonstrate how load shifting can reduce the imbalance, enhancing both grid stability and energy efficiency. Aggregating EV charging can significantly improve renewable energy utilization.}
  \label{fill-fig}
\end{figure}

Figure \ref{fill-fig} illustrates a simplified example of how aggregated EV flexibility could improve grid efficiency. Plot (a) shows an artificial aggregated load profile for public chargers, while plot (b) depicts the corresponding system imbalance reported by the transmission operator, which fluctuates between energy surplus and scarcity. In this case, there is no initial relationship between EV charging and the system imbalance. However, by shifting charging demand from periods of surplus to periods of scarcity (as shown in plot (c)), we can significantly reduce the system imbalance. By redistributing the charging load, we can help absorb excess renewable generation during periods of high solar or wind output, reducing the need for curtailment and improving the utilization of renewable energy sources.

Plot (d) demonstrates that with simple load shifting, the imbalance is reduced by 34.21\%, saving the TSO significant costs in system support fees and enhancing energy efficiency. This kind of aggregated flexibility provides a dual benefit: it helps lower operational costs for grid operators and enhances the sustainability of the energy system by increasing the share of renewable energy used in meeting demand. The cost-effectiveness of this approach makes it an attractive option for both energy providers and consumers, who can benefit from lower energy prices and reduced carbon footprints.

In summary, aggregating EV charging flexibility represents a scalable, cost-effective solution to enhancing energy efficiency and grid stability. By aligning charging with renewable energy availability and reducing peak demand, we can lower both energy costs and emissions. This approach contributes directly to the goals of decarbonizing the grid while supporting the continued growth of renewable energy.

\subsection{Vehicle-to-Grid (V2G) Technology}

The integration of Vehicle-to-Grid (V2G) technology introduces additional flexibility in managing grid imbalances by allowing EVs to discharge energy back to the grid when needed \cite{wevj15110514, HANNAN2022130587}. Unlike traditional charging strategies, V2G enables EVs to serve as dynamic, bidirectional energy resources, providing support during periods of grid stress.

\begin{enumerate}
    \item Grid Support Through Discharging: V2G technology enables EVs to feed electricity back into the grid, offering an additional means of maintaining grid stability. This bidirectional energy transfer allows optimization strategies to account for both charging and discharging periods during times of grid imbalance. For instance, when the grid experiences a shortage of energy supply (i.e., $I_t > 0$), the EV can switch to discharging mode to help balance the system. By redistributing energy in this way, grid imbalances can be reduced, and the total charging requirements for the EV can be minimized through strategic scheduling.

    \item Early-Session Discharging: If an EV is sufficiently charged upon arrival at the charging station, it can begin discharging energy immediately, even in the first hour of the session. This is particularly useful when the grid requires immediate support. However, the optimization framework must carefully adjust the charging schedule based on the vehicle’s initial state of charge (SOC). Discharging at the start of the session can decrease the energy needed to recharge the vehicle later, achieving both immediate grid support and a reduction in overall energy consumption during the session.

    \item Managing Charging and Discharging Under Uncertainty: Since the exact SOC of the vehicle is not always known during the session, the optimization algorithm assumes a conservative, worst-case scenario. This ensures that the total energy delivered to the EV remains within observed limits. For example, if data shows 40 kWh of energy was added to the EV while 10 kWh was discharged during the session, the net energy interaction with the grid would only permit a maximum of 50 kWh to avoid exceeding the vehicle’s capacity. By accounting for discharging in this way, the algorithm ensures the vehicle is not overcharged while staying within operational constraints.
\end{enumerate}
% \newpage
\section{Findings from a Case Study on Real Data from Prague}

In this section, we present the findings from our case study on real-world electric vehicle (EV) charging data from Prague. The study aims to evaluate how flexible EV charging can be optimized to improve energy efficiency and grid stability. By utilizing machine learning (ML) techniques, we provide insights into how EV charging can be shifted to periods of low grid demand, ultimately reducing imbalances and enhancing energy efficiency.

\subsection{Data and Approach}

For our experimentation, we used real-world public EV charging session data from PREdistribuce, a.s.\footnote{Electricity distribution provider in the territory of Prague}, which is the largest operator of public charging infrastructure in Prague. The dataset includes detailed information about charging sessions such as start times, end times, energy transferred, and average power consumption. This granular level of detail allowed us to aggregate charging sessions into hourly loads, which formed the basis of our optimization process. 

To understand the impact of EV charging on grid imbalances, we matched these aggregated loads with system imbalance data provided by ČEPS, the Czech Transmission System Operator (TSO). The system imbalance data represent the difference between electricity supply and demand at any given hour. This information is crucial for identifying periods when demand exceeds supply (requiring upward flexibility) or when there is excess electricity generation (requiring downward flexibility). By matching the charging load with system imbalances, we can explore how EV charging patterns affect overall grid efficiency.

% For this study, we focused on an 8-day period between 25 May 2022 and 2 June 2022, where the data from the public chargers and the system imbalance data had the best overlap. 
For this study, we focused on the periods between 1 June 2022 and 1 July 2022, and 1 June 2024 and 1 July 2024, where the data from the public chargers and the system imbalance data had the best overlap. During this period, we observed typical charging behaviors and imbalances that are representative of real-world conditions. Appendix \ref{ap-charging-data} describes the charging data in more detail, while Appendix \ref{ap-grid-data} provides further insights into the system imbalance data.

This period allowed us to capture both weekday and weekend charging behaviors, which vary significantly due to differences in commuting patterns and public transport usage. Such insights are critical when optimizing charging flexibility, as charging demand fluctuates with human activity. Moreover, by focusing on a real-world dataset, we ensure that our findings are not limited to theoretical models but are instead grounded in the actual behavior of EV users and grid operators.

\subsection{Optimization of Charging}

In the initial analysis, we explored theoretical models to demonstrate the potential of EV charging flexibility in reducing grid imbalances. However, as noted in Section \ref{potential-flex}, such models often rely on unrealistic assumptions. To address this limitation, we performed an optimization based on real-world data to explore more practical and feasible solutions. The optimization sought to shift individual charging sessions within the following realistic constraints \cite{diaz2019optimal}:

\begin{itemize}[nosep]
  \item The start and end of the charging session must remain unchanged, ensuring that the vehicle is fully charged by the time it is needed.
  \item Total power consumption over the session must remain the same, meaning that the optimization focuses solely on the distribution of energy use within the session window.
  \item The vehicle’s state of charge at the end of the optimized session must be the same as in the original data, ensuring that users are not inconvenienced by delayed or incomplete charging.
\end{itemize}

These constraints were designed to ensure that the optimization process was both user-centric and grid-friendly. Importantly, the optimization maintains the user’s expected charging experience while enabling more efficient grid operations by reducing energy use during peak hours.

To perform this optimization, we employed a machine learning model that could redistribute the charging load in a way that minimized the total system imbalance. Specifically, we used a neural network to predict optimal charging times based on historical imbalance data. The model was trained using information about the current and average power imbalance over the past 3 hours, as well as the total session charge consumption and time remaining in the charging session. The goal was to distribute the charging load and enable discharging in such a way that system imbalances were reduced, without negatively impacting the user’s ability to fully charge their vehicle or the vehicle’s operational constraints.

\subsection{Results of the Optimization}

Figures in Appendix \ref{optimization-results} present the results of our optimization efforts. The true observed data are shown in subplots (a) Original observed charging and (b) Original observed imbalance, while the optimized charging sessions are depicted in subplots (c) Optimized charging and (d) Imbalance after optimization. In both the real and optimized scenarios, the charging load pattern is presented in kilowatts (kW), and the system imbalance is shown in megawatts (MW). A summary of the system imbalance areas is also provided in Table \ref{tab:area_system_imbalances}.

\begin{table}[ht]
\caption{System imbalance areas (in MWh) of the Optimal Algorithm during June 2022 and June 2024}
\label{tab:area_system_imbalances}
\centering
\begin{tabular}{|c||c|c|c|c|}
\hline
\textbf{Model} & \textbf{Time Period} & \textbf{Original (MWh)} & \textbf{No V2G (MWh)} & \textbf{V2G (MWh)} \\ \hline
\multirow{2}{*}{Optimal Algorithm} & June 2022 & 366.11 & 365.52 & 365.31 \\ \cline{2-5}
                                   & June 2024 & 402.11 & 400.52 & 400.12 \\ \hline
\end{tabular}
\end{table}

The optimization algorithm shown in this study demonstrated a modest yet statistically significant ability to mitigate grid imbalance, both with and without the integration of V2G technology. While these results highlight the potential of electric vehicles to contribute to grid stability, the study suggests that relying solely on EVs for grid balancing may not be sufficient. A more holistic approach, incorporating diverse flexibility resources such as large-scale batteries, heat pumps, and residential energy storage systems, is likely to yield more substantial and impactful results.

To evaluate the effectiveness of these optimizations, we calculated the percentage of imbalance reduction achieved during June 2022 and June 2024.

For June 2022:
\begin{itemize}
    \item Original imbalance: 366.11 MWh
    \item Optimized imbalance (No V2G): 365.52 MWh
    \item Optimized imbalance (V2G): 365.31 MWh
\end{itemize}

Percentage of imbalance reduced (No V2G):
\[
\left( \frac{366.11 - 365.52}{366.11} \right) \times 100 \approx 0.16\%
\]

Percentage of imbalance reduced (V2G):
\[
\left( \frac{366.11 - 365.31}{366.11} \right) \times 100 \approx 0.22\%
\]

% \newpage

For June 2024:
\begin{itemize}
    \item Original imbalance: 402.11 MWh
    \item Optimized imbalance (Without V2G enabled): 400.52 MWh
    \item Optimized imbalance (With V2G enabled): 400.12 MWh
\end{itemize}

Percentage of imbalance reduced (Without V2G enabled):
\[
\left( \frac{402.11 - 400.52}{402.11} \right) \times 100 \approx 0.39\%
\]

Percentage of imbalance reduced (With V2G enabled):
\[
\left( \frac{402.11 - 400.12}{402.11} \right) \times 100 \approx 0.49\%
\]

In summary, our findings demonstrate that real-world optimization of EV charging can yield substantial benefits in terms of reducing system imbalances and enhancing energy efficiency. The machine learning model proved effective in identifying optimal charging times that both satisfy user needs and minimize grid strain. As EV adoption continues to grow, these findings suggest that integrating smart charging solutions could play a critical role in future energy systems, offering a scalable, cost-effective approach to balancing renewable energy integration with grid reliability.
% \newpage
%\section{Results}

%discussion - it works
% \newpage
\section{Conclusions and Further Work}

This paper has demonstrated how flexible control of public EV charging can serve as an effective method for regulating system imbalances within power grids, playing a key role in the transition towards greener energy systems. The findings are particularly relevant to the Czech TSO's System Services (SyS), which are tasked with maintaining the quality and reliability of electricity supply. By preserving critical parameters such as voltage and frequency through dynamic EV charging, this study highlights the untapped potential of EVs as flexible resources. Their inherent ability to adapt charging times makes them an essential component in balancing services like automatic and manual frequency restoration, mitigating the volatility caused by increasing shares of renewable energy.

While the research is focused on the Czech Republic, the broader implications are clear: the approach and findings are applicable across interconnected European transmission systems, with minor regional adaptations. EVs, due to their increasing presence and their flexible charging potential, represent a scalable, low-cost solution for enhancing grid stability. Moreover, as EV adoption accelerates globally, the methodologies and machine-learning models explored here could help optimize EV-based flexibility services in public charging infrastructure across multiple geographies. These findings not only contribute to improving grid reliability but also help to integrate higher levels of renewable energy sources without compromising system stability.

Beyond balancing services, the contribution of EV flexibility extends to increasing energy efficiency, reducing reliance on fossil-fuel-powered peaking plants, and minimizing the economic costs associated with ancillary services. As EV batteries can act as both storage and flexible load units, they can absorb excess renewable generation during low-demand periods and shift demand away from peak hours, reducing overall energy waste and supporting decarbonization goals. These dual benefits—improved grid stability and enhanced energy efficiency—position EVs as a cornerstone of future energy systems.

Several avenues of further research are planned to expand upon the limitations of this initial study. Firstly, we are expanding the dataset to include more charging and imbalance data beyond 2022 to better generalize the findings across different time periods and grid conditions. By incorporating a more extensive dataset, we aim to capture seasonal variations, different weather patterns, and more diverse user behaviors, which can have significant impacts on grid imbalances and charging flexibility. This will allow for more robust conclusions regarding the scalability and applicability of the proposed solutions.

Secondly, there is potential to refine the technical parameters used in the study. One area of interest is to adjust the system imbalance granularity, which could offer new insights into short-term fluctuations and allow for even more precise optimization of charging patterns. Additionally, we aim to explore the inclusion of vehicle-to-grid (V2G) technology, where EVs not only consume energy but also feed excess power back into the grid during peak times \cite{qin_toward_2023}. V2G integration would significantly enhance the role of EVs as active participants in grid balancing, turning them into bidirectional energy assets rather than mere consumers.

Thirdly, we are examining how the optimization model performs under more realistic conditions where future grid imbalances and charging session lengths are unknown. In real-world applications, grid operators often have to respond to unforeseen events, such as sudden drops in renewable generation or unexpected surges in demand. By incorporating uncertainty into our machine-learning models, we aim to develop more resilient optimization frameworks that can perform well under variable and unpredictable grid conditions.

Finally, this study highlights the need to address the technical and regulatory challenges associated with integrating such a system into existing grid structures. Coordination between TSOs, distribution system operators (DSOs), and EV charger operators will be essential to maximize the benefits of charging flexibility. Developing a regulatory framework that supports dynamic pricing, demand response incentives, and clear protocols for V2G interactions will be critical in unlocking the full potential of this technology. We invite policymakers and the scientific community to engage in further research and discussions on how best to design such frameworks, ensuring that the technical advances made here are aligned with broader energy transition goals.

In conclusion, public EV charging flexibility represents a significant opportunity to enhance grid stability, increase energy efficiency, and support the integration of renewable energy. This study provides a foundational approach to optimizing this flexibility, but the path forward includes addressing data gaps, refining technical models, and resolving regulatory barriers. We encourage continued research and collaboration in this area to ensure that EVs can play a transformative role in the future of energy systems worldwide.

% \newpage

\begin{ack}
The authors would like to thank and acknowledge PREdistribuce, the Prague Distribution System Operator (DSO), for the charging data used in this study, and ČEPS, the Czech Transmission System Operator (TSO) for the open data on system imbalance in Czechia. This work was supported by grants number SGS24/093/OHK5/2T/13 and SGS23/117/OHK5/2T/13 provided by CTU Prague, and number TS01020030 provided by the Technology Agency of the Czech Republic. Authors declare no conflict of interest.

\end{ack}

\medskip

\small

\bibliography{references}
\bibliographystyle{abbrv}

\appendix
\newpage
\section{Data description}

\subsection{Charging data}
\label{ap-charging-data}
\begin{figure}[H]
    \begin{center}
        \includegraphics[width=1\linewidth]{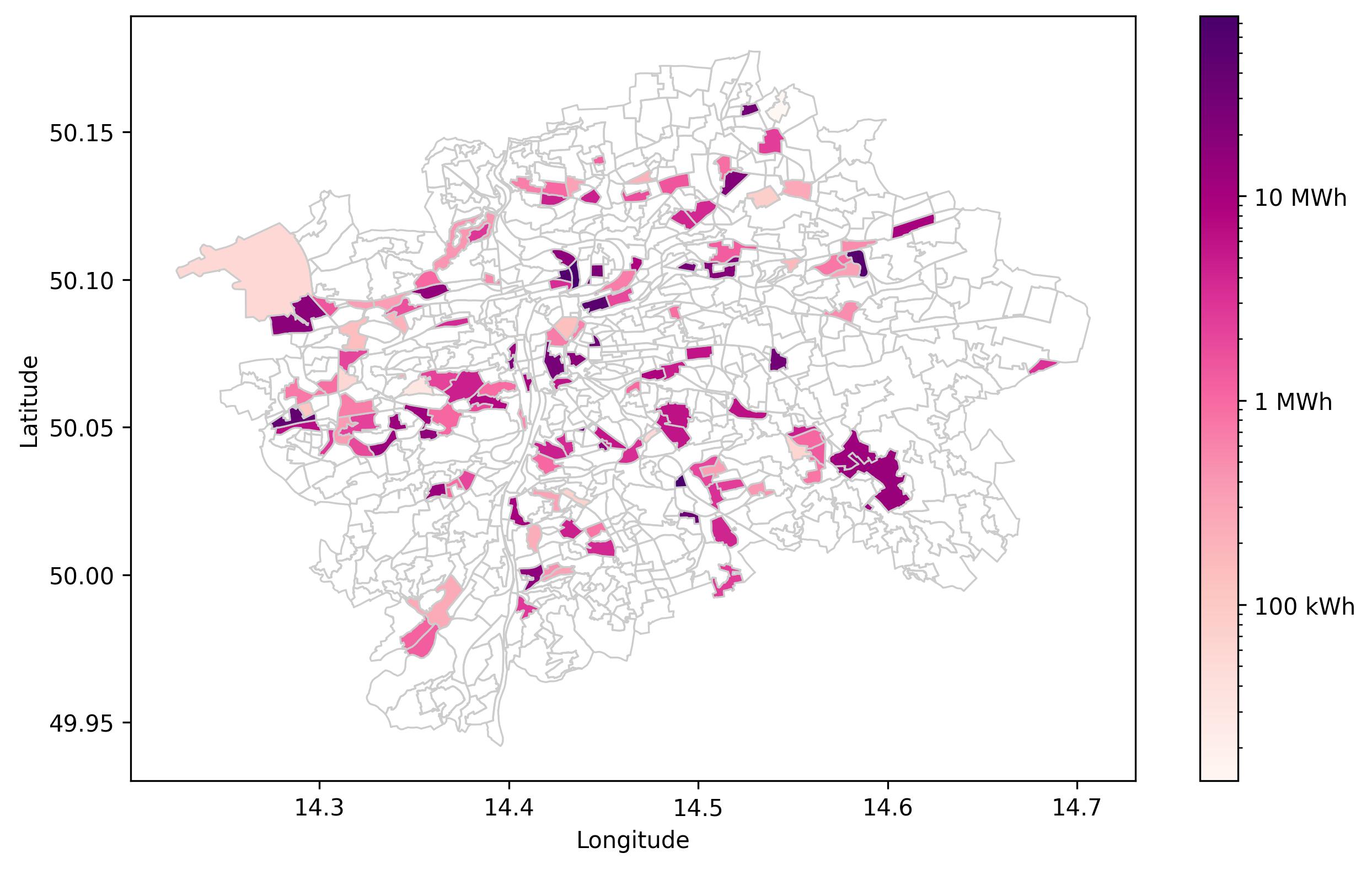}
    \end{center}
                    \vspace{-0.5em}
    \caption{Heatmap of public charging point locations and yearly power consumption per basic administrative unit in Prague based on the available data}
                \vspace{-1em}
    \label{charger-loc-fig}
\end{figure}

The data on individual charging sessions was kindly provided by Prague's largest provider of public charging points, PREdistribuce, which serves between 60-80\% of public chargers 
in Prague.
% \cite{predistribuce_verejne_2024}. 
Figure \ref{charger-loc-fig} shows a visualization of charger locations within the basic administrative unit structure of Prague.

This dataset is currently not publicly available. However, we are preparing to release a sanitized version to the open domain in collaboration with PREdistribuce.

\newpage
\subsection{Grid balance data}
\label{ap-grid-data}
\begin{figure}[H]
    \begin{center}
        \includegraphics[width=1\linewidth]{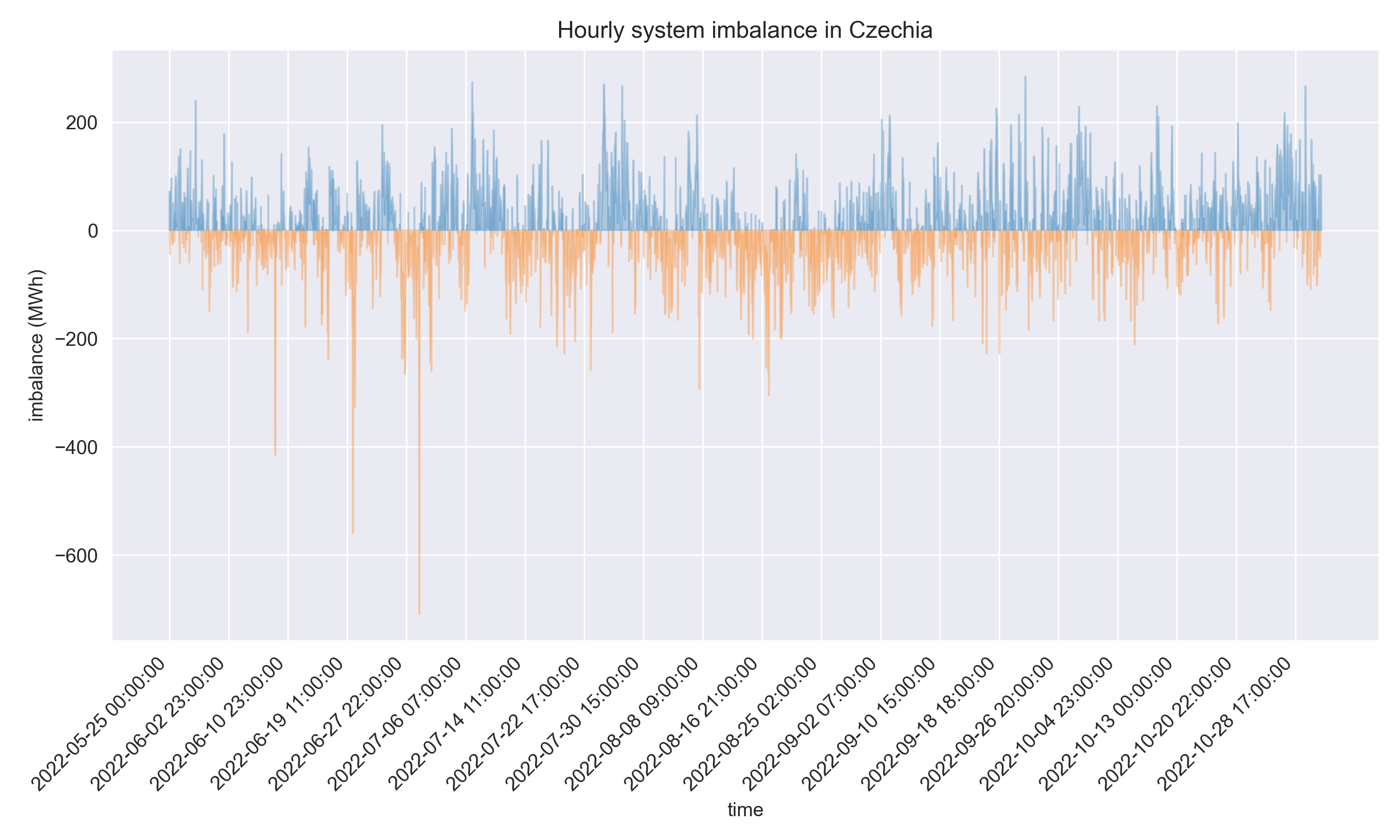}
    \end{center}
                    \vspace{-0.5em}
    \caption{Hourly system imbalance in Czechia between June and October 2022}
                \vspace{-1em}
    \label{system-imbalance}
\end{figure}

Since the charging data we are working with is located in Prague, in this study, we have worked with power grid system imbalance from the Czech TSO, ČEPS. An example 3 month subset of the data can be seen on figure \ref{system-imbalance}. Note that compared to this figure, we have divided the system balance by a factor of 10 in our experiment as the charging loads are from Prague only and Prague represents roughly 10\% of Czech population and power consumption. This is for a slightly easier visual understanding when visualising in one figure along with the aggregated charging loads described in \ref{ap-charging-data}. While not technically accurate as to how system imbalance is geographically distributed in practice, it serves our simplified experiment well as it has minimal impact on results, since the resulting imbalance is still several orders of magnitude larger than the aggregated charging data.

This dataset is freely available as open data updated each day and going several years back, with several periods of aggregation down to individual minutes. It is available on the following url:
\begin{center}
  \url{https://www.ceps.cz/cs/data#AktualniSystemovaOdchylkaCR} 
\end{center}

\newpage
\subsection{Optimization results}
\label{optimization-results}

The plots on the left show the true observed (a) and the optimized (c) charging load pattern in kW. The plots on the right correspond to the real observed (b) and result from the optimization (d) system imbalance attributed to Prague in MW.
% Since the difference would be difficult to perceive in the MW scale, a red line is added to (d) with a secondary kW y-axis on the right, which demonstrates the difference from (b). 
Additionally, the total imbalance in MWh is printed in both (b) and (d) plots.

\subsubsection{Optimization without V2G enabled}

\begin{figure}[H]
    \begin{center}
        \includegraphics[width=1\linewidth]{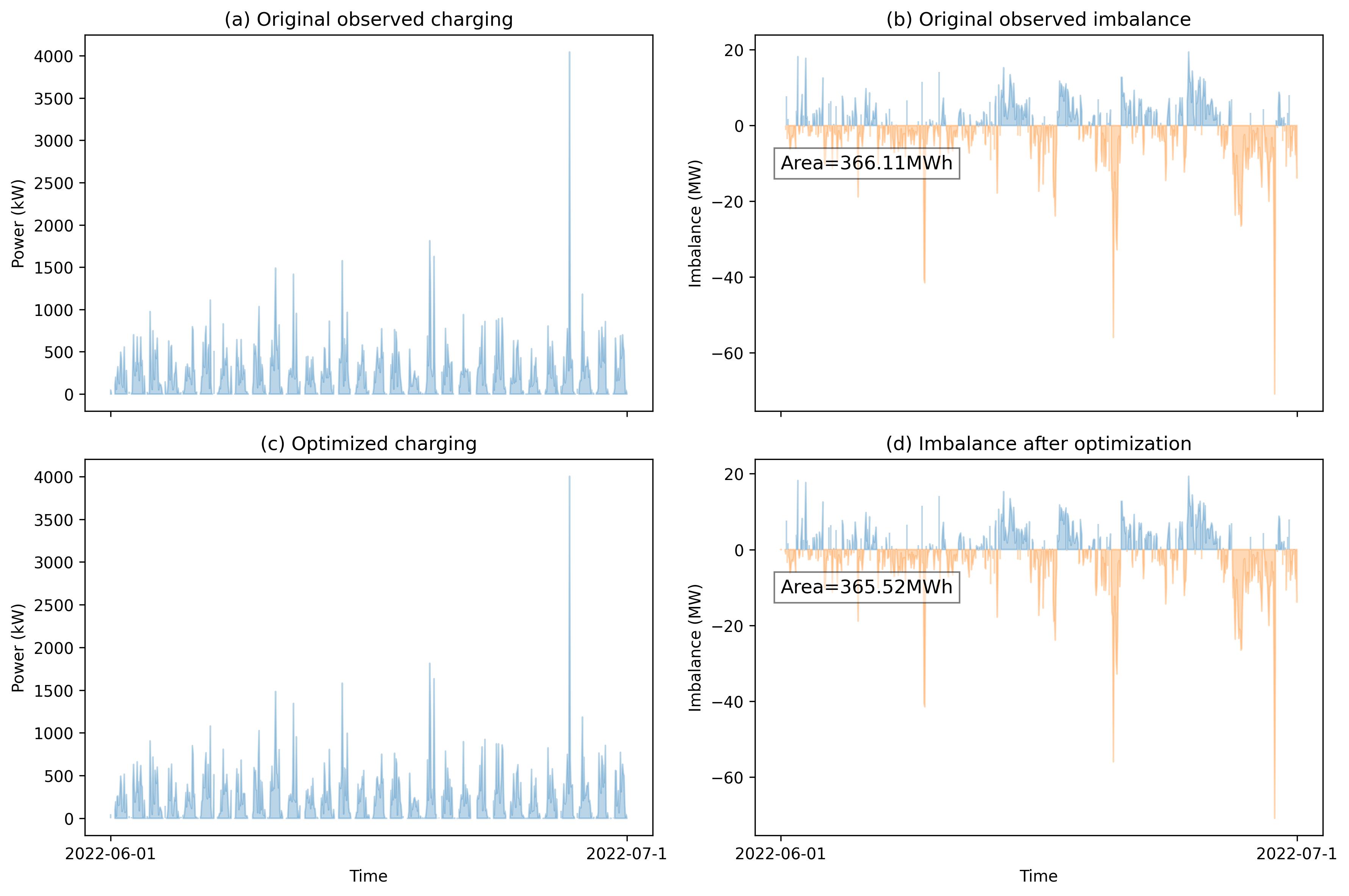}
    \end{center}
                    \vspace{-0.5em}
    \caption{Hourly system imbalance in June 2022 using optimal control under perfect information.}
                \vspace{-1em}
    \label{optimization-june-2022}
\end{figure}

\begin{figure}[H]
    \begin{center}
        \includegraphics[width=1\linewidth]{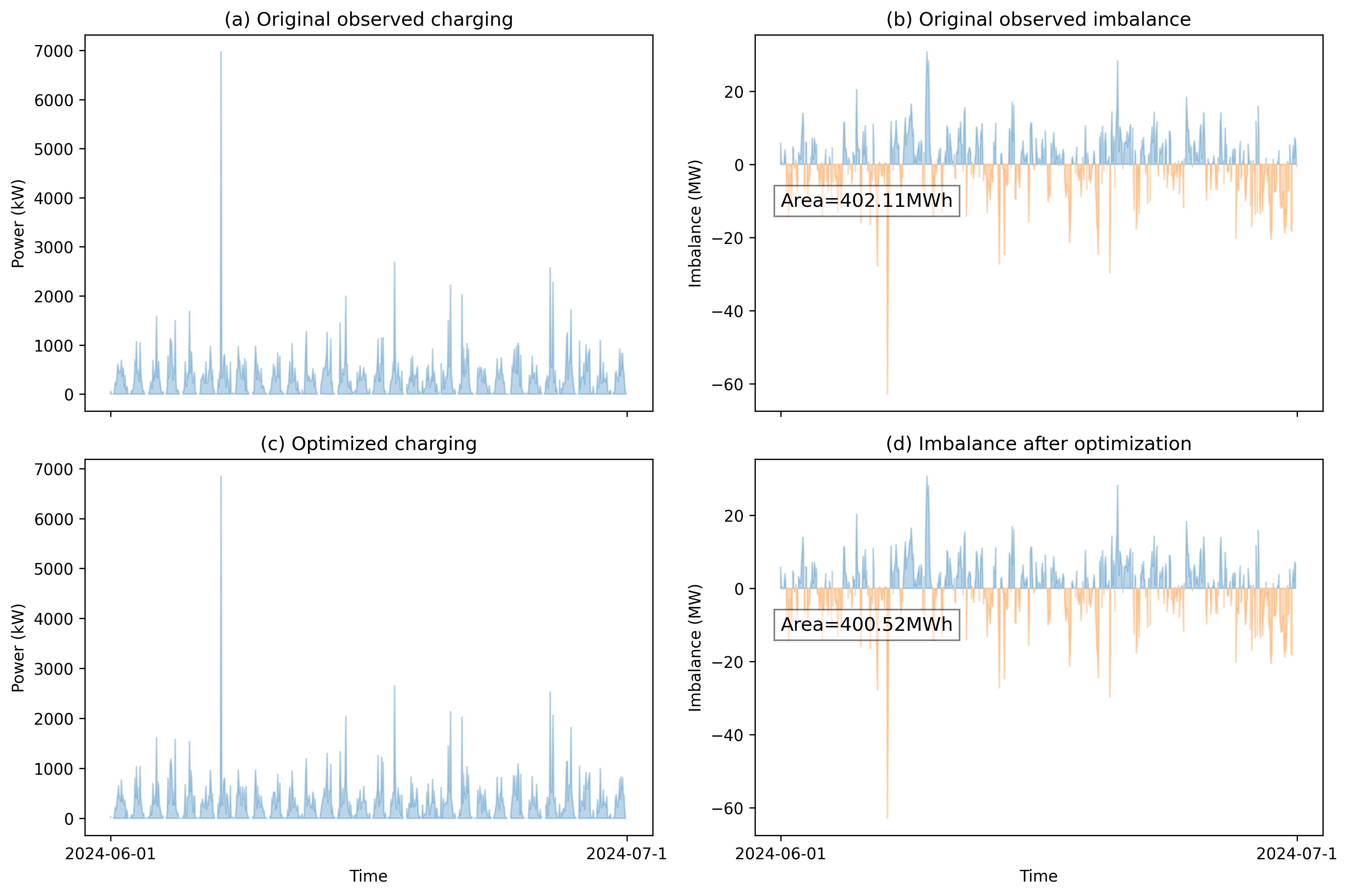}
    \end{center}
                    \vspace{-0.5em}
    \caption{Hourly system imbalance in June 2024 using optimal control under perfect information.}
                \vspace{-1em}
    \label{optimization-june-2024}
\end{figure}

\subsubsection{Optimization with V2G enabled}

\begin{figure}[H]
    \begin{center}
        \includegraphics[width=1\linewidth]{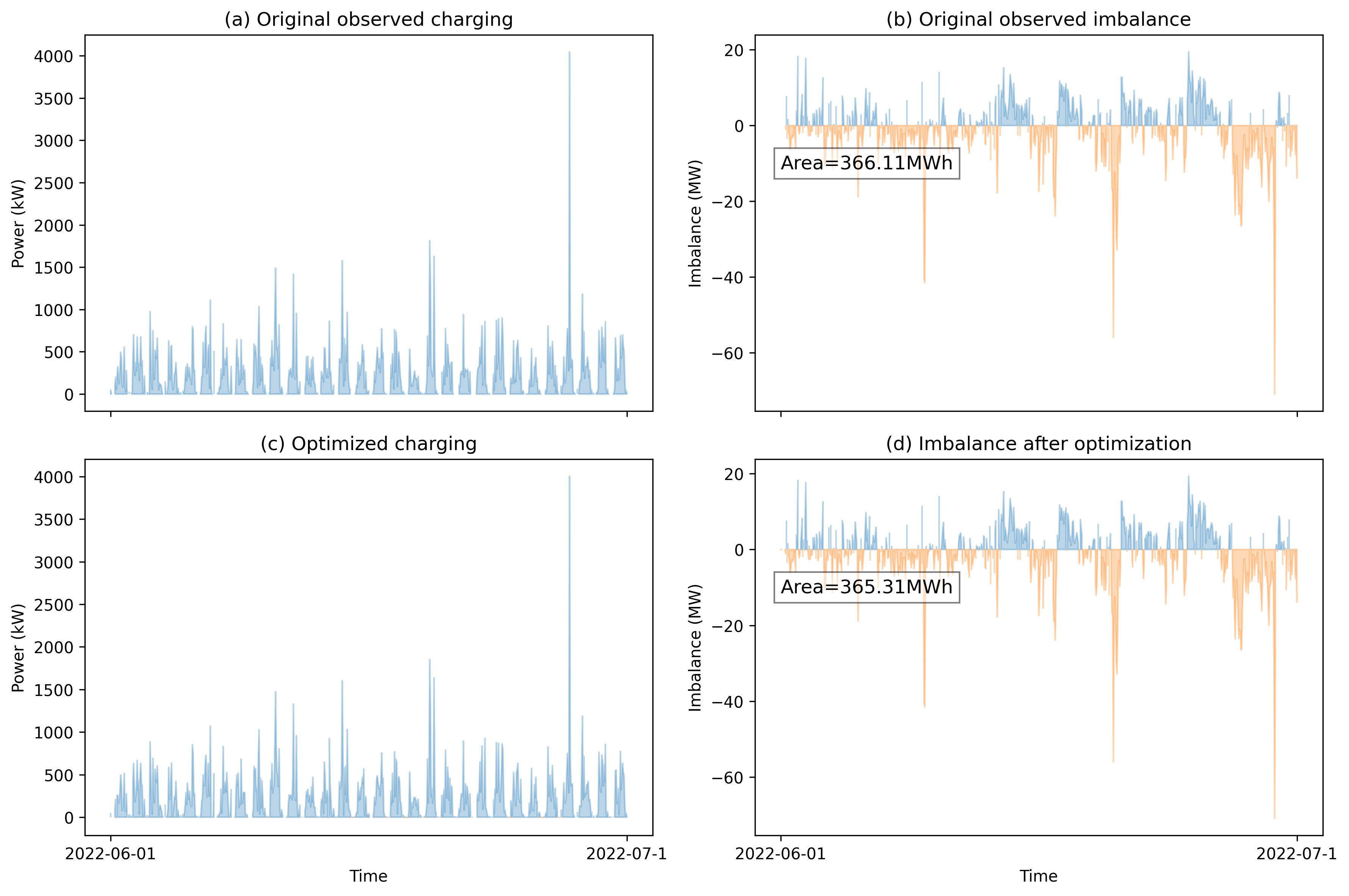}
    \end{center}
                    \vspace{-0.5em}
    \caption{Hourly system imbalance in June 2022 using optimal control under perfect information with V2G enabled.}
                \vspace{-1em}
    \label{optimization-june-2022-v2g}
\end{figure}

\begin{figure}[H]
    \begin{center}
        \includegraphics[width=1\linewidth]{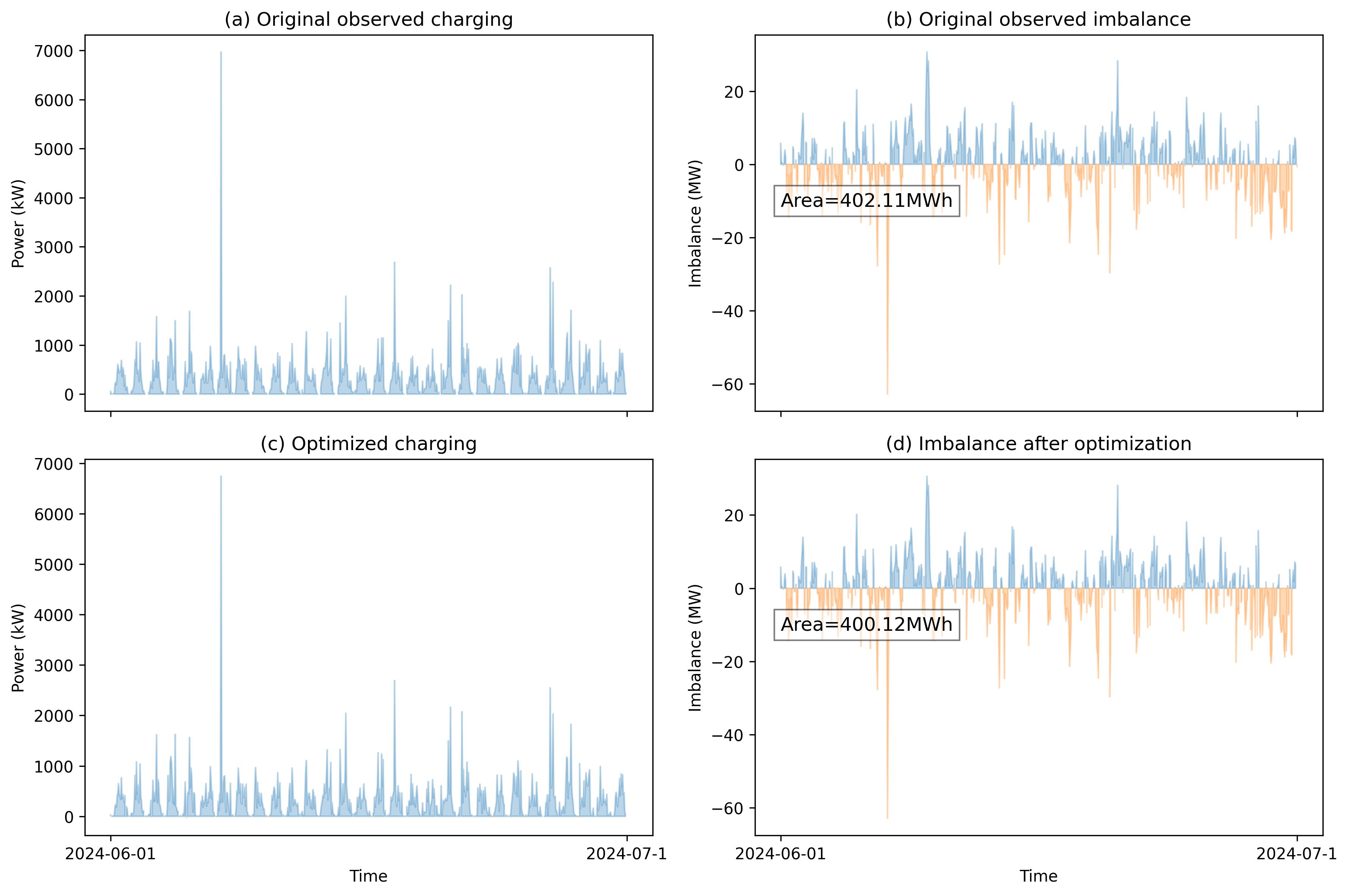}
    \end{center}
                    \vspace{-0.5em}
    \caption{Hourly system imbalance in June 2024 using optimal control under perfect information with V2G enabled.}
                \vspace{-1em}
    \label{optimization-june-2024-v2g}
\end{figure}

\end{document}